\documentclass{article}%

\usepackage{graphicx}
\usepackage{newtxtext,microtype}
\usepackage{newtxmath}
\usepackage{natbib}
\usepackage[hidelinks]{hyperref}
\usepackage{enumerate}
\usepackage{enumitem}
\usepackage{verbatim}
\usepackage{xcolor}

\usepackage[margin=0.85in]{geometry}

\graphicspath{{figures/}}

\def\XXint#1#2#3{{\setbox0=\hbox{$#1{#2#3}{\int}$}
     \vcenter{\hbox{$#2#3$}}\kern-.5\wd0}}

\newcommand*{\e}{\mathrm{e}}
\newcommand*{\ep}{\epsilon}

\newcommand*{\ii}{\mathrm{i}}
\newcommand*{\de}{\mathrm{d}}

\numberwithin{equation}{section}

\begin{document}

\title{Automorphisms of Stokes multipliers in higher-order WKBJ theory}

\author{Josh Shelton$^{1}$$^{\text{\textsection}}$
\footnote{Electronic address: josh.shelton@st-andrews.ac.uk},
\hspace{0.1cm} Samuel Crew$^2$$^{\text{\textsection}}$
\footnote{Electronic address: samuel.c.crew@gmail.com}\hspace{0.1cm} and Christopher J. Lustri$^3$
\footnote{Electronic address: christopher.lustri@sydney.edu.au} {\let\thefootnote\relax\footnote{JS and SC equally contributed to the method developed in this work}}}
\date{$^1$School of Mathematics and Statistics, University of St Andrews, St Andrews, KY16 9SS, UK\\
$^2$Department of Physics, National Tsing Hua University, Hsinchu, Taiwan \\
    $^3$School of Mathematics and Statistics, The University of Sydney, Sydney, NSW, 2006,  Australia\\[2ex]%
}

\maketitle

\abstract{
\noindent We consider the Stokes phenomenon and higher-order Stokes phenomenon (HOSP) of formal asymptotic transseries arising in the WKBJ analysis of linear differential equations and integral problems. We introduce a framework of automorphisms that act on the Stokes constants of the divergent expansion, explained via late-late-term expansions and parametric Alien calculus, to capture this phenomenon. Our method is applied to a paradigmatic example: we obtain the full Stokes line structure and automorphisms for the Swallowtail problem from catastrophe theory, which contains four WKBJ components. We demonstrate that, in a system with four or more WKBJ components, the automorphism associated with the HOSP can itself change value across another higher-order Stokes line, which occurs when different higher-order Stokes lines intersect. We then argue that no additional special behaviour emerges for transseries with five or more WKBJ components.
}

\section{Introduction}
\label{sec:intro}
Stokes’ phenomenon plays a critical role in characterising the asymptotic behaviour of solutions to singularly perturbed problems throughout different sectors of the complex plane. A canonical example occurs in the Airy functions, which transition along the real axis from oscillatory behaviour to uniform exponential decay or growth due to Stokes switching \citep{stokes1847theory,stokes1864discontinuity}. However, this famous example does not display the full range of behaviour that can be seen due to the Stokes phenomenon \citep{berk1982new,howls2004higher,body2005exponential,chapman_2005,shelton2025HOSP,aniceto2024algebraic}.

The Airy functions solve a linear second-order singularly perturbed differential equation whose asymptotic transseries contains exactly two WKBJ contributions. The interaction of these two contributions generates Stokes switching, causing the change in behaviour. In systems with three or more WKBJ components, it is known that the asymptotic structure becomes more complicated, giving rise to the higher-order Stokes phenomenon (HOSP) that causes the effect of Stokes lines to vary throughout the complex plane \citep{howls2004higher,body2005exponential,chapman_2005}.

In this study, we develop a framework to capture automorphisms of the Stokes constants that is used to establish two related results. First, in systems whose transseries contains four or more WKBJ components, the behaviour of the HOSP itself generally changes throughout the complex plane. An example of this phenomenon is also previously noted in \cite{body2005exponential}. Critically, we show that these changes also occur across existing higher-order Stokes lines. Second, no fundamentally new behaviour appears once five or more WKBJ components are present. As a consequence, we show that the full generality of Stokes’ phenomenon is realised in systems with four WKBJ components.
These results are obtained by using the language of Alien calculus, which we use to explain how the higher-order Stokes phenomenon emerges in this framework. 
We interpret the Stokes and higher-order Stokes phenomenon as automorphisms which act on constants associated with the asymptotic transseries. 
The Stokes automorphism acting on the transseries parameters (which multiply each of the WKBJ components) is the Stokes phenomenon, and the size of the resultant jump is proportional to the corresponding Stokes constant.
The higher-order Stokes automorphism, which acts on the Stokes constants (which control the proportion of one WKBJ amplitude function present in the late-terms of another) is the higher-order Stokes phenomenon, and the size of the resultant jump is shown to be proportional to the product of other Stokes constants.

\subsection{The Regular and Higher-Order Stokes Phenomenon}
The Stokes phenomenon can be understood by studying the full asymptotic expansion of some function arising as a solution to an ODE (its \emph{asymptotic transseries}) in the complex plane. Consider a complex-valued function $\psi$ in terms of a complex variable $z$, and which involves some parameter $\epsilon$.
In the limit as $\epsilon \to 0$, the solution takes the form of an asymptotic transseries
\begin{equation}\label{eq:transseries0}
    \psi(z,\ep) \sim \sum_{i=1}^{N}  \epsilon^{-\alpha_i} \sigma_i \e^{-\chi_{i}(z)/\ep} \Big[\psi^{(i)}_0(z)+\ep \psi^{(i)}_1(z)+\cdots \Big],
\end{equation}
where $\chi_i$ and $\psi_i$ are complex functions, which consists of every possible order of the asymptotic expansion. Each distinct $\chi_i$ represents a different WKBJ component. The transseries parameters $\sigma_{i}$ are piecewise-constant, but change discontinuously in different sectors of the complex plane bounded by Stokes lines. This discontinuous variation produces the rapid change in exponential contributions associated with the Stokes phenomenon. A more detailed discussion of resurgent transseries and their form can be found in \cite{aniceto2015nonperturbative}.

The Stokes phenomenon appears in asymptotic transseries with at least two different WKBJ components, such as those which arise from second-order singularly-perturbed linear differential equations. This context is well-understood, and related studies that build on this understanding of the Stokes phenomenon include the study of nonlinear differential equations solved via divergent expansions and matched asymptotics \citep{chapman1998exponential} or WKBJ theory \citep{kawai2005algebraic}, saddle-point integrals solved by hyperasymptotics \citep{berry1991hyperasymptotics}, and Gevrey-divergent transseries expansions via \'Ecalle's theory of Alien calculus \citep{aniceto2019primer} which builds upon the fundamental divergence of asymptotic expansions first established by \cite{dingle1973asymptotic}.

The rapid switching associated with the Stokes phenomenon was shown by \cite{berry1989uniform} to occur in a smooth manner across the Stokes line, and rigorous error bounds may be obtained for asymptotic transseries that display Stokes phenomenon by studying the hyperterminant functions that emerge in the theory of hyperasymptotics \citep{daalhuis1996hyperterminants,daalhuis1998hyperterminants}. The asymptotic form \eqref{eq:transseries0} possesses additional subtleties due to the behaviour of $\sigma_i$. The change in value of $\sigma_i$ across Stokes lines is localised to a boundary layer of diminishing width under the asymptotic limit; for most considerations, only the total change in value across the Stokes line is required. It is therefore typical to regard the transseries parameters as piecewise constant functions whose value differs either side of the Stokes line, for which the Stokes phenomenon can be regarded as a linear automorphism acting on the set of transseries parameters.

These studies largely discuss regular Stokes switching, associated with the interaction between two WKBJ components. However, more complicated behaviour emerges as the number of interacting components increases. The \emph{higher-order} Stokes phenomenon (HOSP), identified by \cite{berk1982new}, manifests in asymptotic transseries with three or more WKBJ components, such as those arising from third-order homogeneous linear differential equations. HOSP affects the activity of regular Stokes lines depending on the location in the complex plane, altering the Stokes switching seen in the asymptotic solution. Further studied by \cite{howls2004higher,body2005exponential,chapman_2005}, the HOSP resolves contradictions associated with intersecting Stokes lines and can result in inactive or new Stokes lines emerging.
More recent studies have focused on deriving the HOSP from the late-terms of the asymptotic expansion \citep{shelton2025HOSP}, and investigating how it influences the smooth behaviour of the original asymptotic expansion across Stokes lines \citep{howls2025smoothing}.

In this paper we demonstrate that the activity of higher-order Stokes lines can vary throughout the complex plane, which occurs for asymptotic transseries with four or more WKBJ components (such as from linear differential equations of order four or higher).
Four or more WKBJ components are required to observe this asymptotic effect, because this phenomenon requires the crossing of two different higher-order Stokes lines to occur, and takes place across other higher-order Stokes lines already present in the expansion. 

The fact that this change occurs across other higher-order Stokes lines leads to an additional consequence, which is that no further new phenomenon appears in transseries expansions with five or more WKBJ components. We will show this by studying the late-late terms of the asymptotic expansion (subdominant components of the late-terms), an interpretation of which is also given in the language of \'Ecalle's Alien calculus \citep{ecalle1981iii}.  

\subsection{The swallowtail Integral}
To demonstrate the necessity of considering higher-order automorphisms in determining Stokes multipliers for a divergent series expansion, we study the Stokes structure for the swallowtail problem from catastrophe theory, 
\begin{equation}
\label{eq:swallowint0}
\Psi(x_1,x_2,x_3) = \int_{-\infty}^{\infty} \exp{\big(\mathrm{i}[t^5+x_3 t^3+x_2 t^2+x_1 t]\big)}\mathrm{d}t,
\end{equation}
where $x_1$, $x_2$, and $x_3$ are real-valued parameters. This integral is a member of the hierarchy of special functions that also contains the Airy and Pearcey functions, and satisfies a fourth-order differential equation \eqref{eq:swallowdiff} that will be introduced in Section \ref{sec:swallowint}. 

The behaviour of the swallowtail integral has previously been studied in certain parameter regimes via both asymptotic and numerical methods.
Note that when studying the swallowtail problem in \S\,\ref{sec:swallowtail} we perform a far-field rescaling of each parameter $x_1$, $x_2$, and $x_3$ in \eqref{eq:swallowint0} to realise the full generality of Stokes phenomenon that emerges due to the inclusion of four saddle points in the integral.
The change in value of Stokes constants that we report has not been detected in the swallowtail problem in any previous work.

The swallowtail integral was evaluated numerically by \cite{connor1983differential,connor1984uniform} for finite values of $x_1$, $x_2$, and $x_3$.
The swallowtail caustic was studied by \cite{wright1981singularities}, and asymptotic solutions were obtained by \cite{kaminski1992asymptotics} in the far field at specific parameter values corresponding to the cusps of this caustic.
More recently, both \cite{ferreira2018asymptotic} and \cite{bennett2018globally} have studied the swallowtail integral with specific parameter values that resulted in three saddle points. A hyperasymptotic method was developed by \cite{bennett2018globally} and used to obtain values for the Stokes constants (often also called multipliers) in this formulation.
In our formulation we will study the case for large $|x_1|$, $x_2=O(x_1^{3/4})$, and $x_3=O(x_1^{1/2})$, which arises naturally when considering the full generality of behaviour produced by the integral with four distinct saddle points. Later, in figure~\ref{fig:abvalues} we show the parameter values for the Pearcey integral at which saddle points coalesce---the cusp of the caustic then corresponds to the origin in this figure.

As the integral is generated using a fourth-order differential equation, it possesses four WKBJ components, and hence contains numerous higher-order Stokes lines. Our analysis will reveal that when different certain higher-order Stokes lines intersect, a variant of the Stokes phenomenon occurs in the higher-order Stokes switching itself. This results in certain higher-order Stokes lines being inactive, and is necessary to take this behaviour into account in order to correctly determine the full Stokes structure of the asymptotic expansion.

As this example contains four WKBJ components, it captures all possible variations of the Stokes phenomenon that can appear in an asymptotic transseries. Equations arising from higher-order differential equations, with five or more WKBJ components, are shown in \S\,\ref{sec:fiveormore} to not introduce any new variants of Stokes switching behaviour.

\subsection{Paper Structure}

We begin in \S\,\ref{sec:method} by deriving the Stokes phenomenon and the higher-order Stokes phenomenon for a Gevrey-divergent transseries expansion.
These phenomena are captured by considering automorphisms that act on both the transseries parameters and the Stokes constants of the expansion.
These automorphisms are developed by using the parametric Bridge equation from Ecalle's theory of alien calculus.
In \S\,\ref{sec:swallowtail} we consider the transseries expansion for the fourth-order swallowtail ODE. The components of the transseries are derived, namely the singulant functions and amplitude functions. We also specify values for the Stokes constants in a certain sector of the complex plane for a given, fixed, choice of boundary conditions.
In \S\,\ref{sec:stokesstructure} we calculate the full Stokes line structure for the swallowtail example by first obtaining the Stokes constants throughout the complex plane by applying our higher-order Stokes automorphisms, and then we obtain the transseries parameters throughout the complex plane from the regular Stokes automorphisms.
We conclude in \S\,\ref{sec:conc}.

\section{Transseries in singular perturbation theory}
\label{sec:method}

Asymptotic solutions to singular perturbative problems typically take the form of an asymptotic transseries, shown previously in \eqref{eq:transseries0},
which consists of every possible order of the asymptotic expansion and holds under the limit of $\ep \to 0$. Such asymptotic transseries frequently arise in the study of differential equations, integral equations, and observables in quantum field theory. Following \cite{dingle1973asymptotic} we refer to the functions $\chi_i(z)$ as singulants, and these will determine where in the complex plane Stokes phenomenon occurs.

In the context where $\psi(z;\epsilon)$ is the solution to some differential equation, the functions $\psi_{n}^{(i)}(z)$ and $\chi_{i}(z)$ can be determined from the governing equation. The transseries parameters $\sigma_i(\epsilon)$ are then set in a certain sector of $z \in \mathbb{C}$ by far-field or behavioural conditions of the problem. The number of linearly-independent WKBJ modes, $N$, in the transseries \eqref{eq:transseries0} is problem specific. For linear ODEs $N$ is typically finite, and equal to the order of the differential equation unless external forcing is applied. For nonlinear ODEs $N = \infty$ in general 
\citep{olde2022exponentially}.

Depending on the value of the singulants $\chi_i(z)$, different modes can exchange dominance with one another as $z$ changes value, which occurs across anti-Stokes lines.
Even in contexts where no significant exchange of dominance occurs to the asymptotic behaviour of $\psi(z;\epsilon)$ across anti-Stokes lines, it is still often necessary to have knowledge of subdominant components.
This occurs for instance when enforcing quantization conditions on subdominant components of the expansion \citep{sueishi2020exact}.
Explaining the behaviour of a complex function in terms of its transseries therefore requires not only tracking exponentially subdominant terms, but also the Stokes phenomenon induced by their interactions.

\subsection{The Stokes phenomenon}
\label{sec:methodlevel1}

The Stokes phenomenon results in distinct asymptotic approximations for the solution $\psi(z;\epsilon)$ holding in different sectors of $z \in \mathbb{C}$. To determine the full transseries expression for a complex function, it is necessary to determine the asymptotic behaviour in each sector, and to identify how to map between the asymptotic approximations as the sector boundaries (ie. Stokes lines and anti-Stokes lines) are crossed. In this section, we describe how to locate Stokes lines, and how to determine the change in $\sigma_i$ as Stokes lines are crossed.

We now consider the Stokes phenomenon acting on the $j$th exponential, $\mathrm{e}^{-\chi_j/\ep}$, in transseries \eqref{eq:transseries0}, induced by the $i$th exponential, $\mathrm{e}^{-\chi_i/\ep}$.
The Stokes phenomenon is a consequence of divergence in the amplitude expansion $\psi^{(i)}(z)\sim \psi_0^{(i)}(z)+\epsilon \psi_1^{(i)}(z)+\cdots$.
For a Gevrey-$1$ asymptotic series that arises as the solution to a linear ODE, the $n$th term of the amplitude expansion diverges as $n \to \infty$ in the factorial-over-power manner of
\begin{equation}\label{eq:lateterms}
    \psi_n^{(i)}(z) \sim \sum_{j} \left[\sum_{p=0}^{\infty}\tilde{\psi}_p^{(j)}(z)\frac{\Gamma\big(n+\alpha_j-\alpha_i-p\big)}{\big[\chi_j(z)-\chi_i(z)\big]^{n+\alpha_j-\alpha_i-p}}\right].
\end{equation}
In \eqref{eq:lateterms}, the summation over $j$ captures different divergent contributions to the late-term asymptotics of $\psi_n^{(i)}(z)$, and lower-order components (in $n$) of each of these are captured by the summation over $p$.

The functions $\tilde{\psi}_p^{(j)}(z)$ appearing in the late-term expansion \eqref{eq:lateterms} may be grouped into an algebraic $\epsilon$ expansion by introducing the Alien derivative, $ \Delta_{\chi_j-\chi_i}$, which acts on the original amplitude expansion $\psi^{(i)}(z)\sim \psi_0^{(i)}(z)+\epsilon \psi_1^{(i)}(z)+\cdots$.
This yields
\begin{equation}\label{eq:originalalien}
    \Delta_{\chi_j-\chi_i}\psi^{(i)}(z) = \sum_{p=0}^{\infty}\epsilon^p \tilde{\psi}^{(j)}_p(z),
\end{equation}
which defines the so-called parametric Alien derivative.\footnote{See for example \cite{aniceto2019primer} and \cite{dorigoni2019introduction} for reviews of Alien calculus in the non-parametric setting.}
It is a linear map which takes a formal Gevrey-1 asymptotic series, $\psi^{(i)}(z)$, to another formal Gevrey-1 asymptotic series, $\tilde{\psi}^{(j)}(z)$.
The Alien derivative satisfies the Leibniz identity for differentiation and, by defining $\dot{\Delta}_{\chi} = \mathrm{e}^{-\chi(z)/\epsilon} \Delta_{\chi}$, the property $ \Dot{\Delta}_{\chi}(\frac{\mathrm{d}F}{\mathrm{d}z}) = \frac{\mathrm{d}}{\mathrm{d}z} ( \dot{\Delta}_{\chi} F )$ holds.
This property implies that when the transseries \eqref{eq:transseries0} arises as a formal solution to a linear differential equation, the Alien derivative of $\psi(z;\epsilon)$ with respect to $\chi_j-\chi_i$ is proportional to the $j$th component of transseries \eqref{eq:transseries0}. Further details and properties of the Alien calculus are provided in Appendix~\ref{sec:alienapp}. This property is often expressed compactly as the (second) bridge equation
\begin{equation}\label{eq:bridge2eq}
\dot{\Delta}_{\chi_j-\chi_i} \psi =\frac{S_{ij}\sigma_i}{2 \pi \mathrm{i}}\frac{\partial \psi}{\partial \sigma_j}, 
\end{equation}
introduced by \cite{ecalle1981iii}. The constants of proportionality $S_{ij}$ are known as the Stokes constants and encode the hard analytic data of the problem. Note that the factor of $(2 \pi \mathrm{i} )^{-1}$ included in \eqref{eq:bridge2eq} ensures that the Stokes constants match with the definitions from \cite{berry1991hyperasymptotics,bennett2018globally}. 
We treat these constants as piecewise constant in the complex-$z$ plane. A consequence of Ecalle's second bridge equation is that the functions $\tilde{\psi}_p^{(j)}(z)$ appearing in the late-term expansion \eqref{eq:lateterms} are equal to the product of the Stokes constant $S_{ij}$ and the amplitude functions $\psi^{(j)}_p(z)$ from the original transseries expansion \eqref{eq:transseries0}, which yields
\begin{equation}\label{eq:tilderel}
 \tilde{\psi}_p^{(j)}(z)=\frac{S_{ij}}{2 \pi \mathrm{i} }  \psi^{(j)}_p(z).
\end{equation}

Each component, indexed by $j$, of the divergence \eqref{eq:lateterms} of $\psi_n^{(i)}(z)$ induces Stokes phenomenon in the asymptotic approximation of the amplitude function $\psi^{(i)}(z)$.
Stokes lines are located wherever the quantity $\chi_j(z)-\chi_i(z)$ is real and positive, and across these a component of the form
$$S_{ij} \epsilon^{\alpha_i-\alpha_j} \psi^{(j)}(z)\mathrm{e}^{-[\chi_j(z)-\chi_i(z)]/\epsilon}$$ 
is gained in the asymptotic expansion of $\psi^{(i)}(z)$.
Upon multiplying this term by $\epsilon^{-\alpha_i} \sigma_i \mathrm{e}^{-\chi_i(z)}$, we see that the overall contribution to the asymptotics of $\psi(z;\epsilon)$ is 
$$S_{ij} \epsilon^{-\alpha_j} \sigma_i \psi^{(j)}(z)\mathrm{e}^{-\chi_j(z)},$$ 
which is a constant multiple of the $j$th term of transseries \eqref{eq:transseries0}.
Thus, the $i$th term of the transseries generates a Stokes phenomenon that acts on the $j$th term, and results in the transseries parameter $\sigma_j$ changing value to $\sigma_{j}+ S_{ij}\sigma_{i}$.
This change in value of the $j$th transseries parameter $\sigma_j$ may be regarded as an automorphism, which we denote by $\mathcal{S}_{i>j}$.
We denote the Stokes line across which this occurs by $l_{i>j}$.
Thus, the Stokes phenomenon consists of a linear automorphism of the transseries coefficients $\sigma_i$ which occurs across the Stokes lines, both of which are defined by
\begin{subequations}
\begin{align}
\label{eq:level1line}
\text{(Stokes line $l_{i>j} \subset \mathbb{C}_z$)}: &\quad  \text{Im}[\chi_j(z)-\chi_i(z)]=0 \quad \text{and} \quad \text{Re}[\chi_j(z)-\chi_i(z)]\geq 0,\\
\label{eq:level1auto}
\text{(Stokes automorphism $\mathcal{S}_{i>j}$)}: &\quad  \sigma_{j} \mapsto \sigma_{j}+ S_{ij}\sigma_{i}.
\end{align}
\end{subequations}
Specifically, the automorphism \eqref{eq:level1auto} acts linearly on the vector of transseries parameters $(\sigma_1,\ldots,\sigma_N)^T$ as the Stokes line $l_{i>j}$ is crossed in the direction $\text{Im}[\chi_j-\chi_i]<0$ to $\text{Im}[\chi_j-\chi_i]>0$. The Stokes automorphism $\mathcal{S}_{i>j}$ may thus be expressed as a matrix with a single off-diagonal entry $S_{ij}$ in row $j$ and column $i$, as given by
\begin{equation}
    (\mathcal{S}_{i>j})_{kl} = \delta_{kl} + S_{ij}\,\delta_{kj}\,\delta_{li}.
\end{equation}

The Stokes lines $l_{i>j}$ are topological half (real) lines emanating from turning points $\tilde{z}$ where $\chi_i(\tilde{z}) = \chi_{j}(\tilde{z})$.
The set of all Stokes lines $\bigcup_{(i,j)} l_{i>j}$ define a Stokes graph $\mathcal{G}_{S}$ as a subset of (a possible multi-sheeted cover of) $\mathbb{C}_z$.

The Stokes phenomenon is now well established in resurgence and WKBJ theory. The asymptotics in \eqref{eq:lateterms}, which describe the large-$n$ behaviour of the coefficients of a resurgent complex function, display Gevrey-$1$ divergence as the expression is bounded by factorial-over-power growth. 
Many equivalent methods exist for deriving this behaviour, including Borel resummation \citep{aniceto2019primer}, saddle-point analysis of integral formulations \citep{berry1991hyperasymptotics,bennett2018globally}, and matched asymptotic expansions \citep{olde1995stokes}.

The value of the Stokes constants $S_{ij}$ is problem dependent (in the language of \'Ecalle, $S_{ij}$ are analytic invariants of the underlying ODE). This is in contrast to $\sigma_i$, which are typically all determined from a behavioural condition in a certain sector of $z \in \mathbb{C}$.
The values of $S_{ij}$ are given formally by the formula
\begin{equation}\label{eq:Stokesconstantform}
S_{ij}=\frac{2 \pi \mathrm{i}}{\psi_0^{(j)}(z)}\lim_{n \to \infty} \left[\frac{\psi_n^{(i)}(z)[\chi_j(z)-\chi_i(z)]^{n+\alpha_{j}-\alpha_{i}}}{\Gamma(n+\alpha_j-\alpha_i)} \right],
\end{equation}
which holds near to the turning point $\tilde{z} \in \mathbb{C}$ at which $\chi_j(\tilde{z})-\chi_i(\tilde{z})=0$.
Since the governing problem can have multiple turning points, it is typical for the values of $S_{ij}$ to be known in different sectors of $z \in \mathbb{C}$. 
For a problem which only consists of the regular Stokes phenomenon, the value of $S_{ij}$ remains unchanged throughout the entire complex plane.
However, for problems with additional turning points it is possible for the value of $S_{ij}$ to change across higher-order Stokes lines, which subsequently changes the value of the multiplier in the Stokes automorphism \eqref{eq:level1auto}. We say that a Stokes line $l_{i>j}$ if \textit{inactive} is $S_{ij}=0$ in some region, and \textit{active} otherwise.

We note that it is possible to determine the value (with rigorous error bounds) of every Stokes constant $S_{ij}$ at a fixed location in the complex plane by the method of \cite{daalhuis1998hyperasymptotic,bennett2018globally}, which determines the Stokes constants to arbitrary precision using hyperasymptotic expansions.

\subsubsection{Adjacency graph}
The (in-)activity of Stokes lines may be conveniently encoded in an adjacency graph \citep{berry1991hyperasymptotics}.
We consider the singulants $\{\chi_1,\ldots,\chi_N\}$ as vertices of a graph and draw an edge if and only if $S_{ij}$ is non-zero. Equivalently, a Stokes line is active if and only if there is a corresponding edge in the graph. The adjacency graph is local in the sense that the edges of the graph vary with $z \in \mathbb{C}_z$. In particular, there is a distinct adjacency graph for each higher-order Stokes region, which we show in detail in the following section.

\subsection{The higher-order Stokes phenomenon}
\label{sec:methodlevel2}
When there are three or more WKBJ components present in \eqref{eq:transseries0}, the asymptotic expansion may exhibit behaviour known as the higher-order Stokes phenomenon, which may be understood as a Stokes phenomenon in the $n \to \infty$ approximation of the late-terms \eqref{eq:lateterms} \cite{shelton2025HOSP,howls2025smoothing}. This results in a change in value of the Stokes constant $S_{ij}$ across curves known as higher-order Stokes lines. Throughout this section, we will describe the interactions caused by three distinct exponentials: $\mathrm{e}^{-\chi_i/\ep}$, $\mathrm{e}^{-\chi_j/\ep}$, and $\mathrm{e}^{-\chi_k/\ep}$ with $i \neq j \neq k$.

The higher-order Stokes phenomenon can be captured through the consideration of lower-order components of the late-term divergence of $\psi_n^{(i)}(z)$ as $n \to \infty$, which was given previously in \eqref{eq:lateterms}.
By substituting for the solution $\tilde{\psi}_p^{(j)}=S_{ij} (2 \pi \mathrm{i})^{-1} \psi_p^{(j)}$ from \eqref{eq:tilderel}, the late-term asymptotics of $\psi_n^{(i)}$ are given by
\begin{equation}\label{eq:latetermsall}
    \psi_n^{(i)}(z) \sim \sum_{j} \frac{S_{ij}}{2 \pi \mathrm{i}} \Bigg[\sum_{p=0}^{\infty}\psi_p^{(j)}(z)\frac{\Gamma (n+\alpha_j-\alpha_i-p )}{\big[\chi_j(z)-\chi_i(z)\big]^{n+\alpha_j-\alpha_i-p}}\Bigg]
\end{equation}
as $n \to \infty$. 
Since $\psi_p^{(j)}(z)$ diverges as $p \to \infty$, the $n \to \infty$ asymptotic expansion of the late-terms \eqref{eq:latetermsall} itself forms a divergent asymptotic series for large $p$.
The divergence of $\psi_p^{(j)}(z)$ is analogous to that for $\psi_n^{(i)}(z)$ in \eqref{eq:lateterms}, and takes the form
\begin{equation}\label{eq:latetermsp}
    \psi_p^{(j)}(z) \sim \sum_{k} \frac{S_{jk}}{2 \pi \mathrm{i}} \psi_0^{(k)}(z)\frac{\Gamma ( p+\alpha_k-\alpha_j )}{\big[\chi_k(z)-\chi_j(z)\big]^{p+\alpha_k-\alpha_j}}
\end{equation}
as $p \to \infty$. In \eqref{eq:latetermsp} we have written for notational convenience only the dominant divergent component for each term in the $k$ summation.
This reexpansion of the late-terms, \textit{i.e.} the late-late-terms of the expansion, may be expressed in terms of the Alien derivative \eqref{eq:originalalien} as
\begin{equation}
    \Delta_{\chi_k-\chi_j}\Delta_{\chi_j-\chi_i}\psi^{(i)}(z) = -\frac{S_{ij}S_{jk}}{(2 \pi)^2}\psi^{(k)}(z),
\end{equation}
where we have applied the bridge equation twice (which notably relies on the assumption that the transseries arises as a formal transseries solution to an ODE).
Note that crucially $\Delta_{\chi_k -\chi_i}\psi^{(i)}(z)$ need not be equal to $\Delta_{\chi_k-\chi_j}\Delta_{\chi_j-\chi_i}\psi^{(i)}(z)$.

The higher-order Stokes phenomenon induced by divergence in the full late-term expansion \eqref{eq:latetermsall} was derived by \cite{shelton2025HOSP} by resumming the remainder of an optimally truncated $p$ summation in \eqref{eq:latetermsall}, and this result has also been generalised to the resultant smooth behaviour on the $\epsilon$ expansion by \cite{howls2025smoothing}. 
The result is that the Stokes constants change value across a higher-order Stokes line, which occurs whenever the quantity $(\chi_k-\chi_j)/(\chi_j-\chi_i)$ is real and positive.
We thus define the higher-order Stokes line $h_{i>k;j}$ and higher-order Stokes automorphism $\mathcal{T}_{i>k;j}$ by
\begin{subequations}
\label{eq:HOSA+P}
\begin{align}
\text{(higher-order Stokes line $h_{i>k;j}  \subset \mathbb{C}_z$)}: \quad  \text{Im}\left[\frac{\chi_k-\chi_j}{\chi_{j}-\chi_{i}}\right]=0 \quad \text{and} \quad \text{Re}\left[\frac{\chi_k-\chi_j}{\chi_{j}-\chi_{i}}\right]\geq 0,\\
\label{eq:HOSA+P2}
\text{(higher-order Stokes automorphism $\mathcal{T}_{i>k;j}$)}: \quad  S_{ik} \mapsto S_{ik}+ S_{ij} S_{jk},
\end{align}
\end{subequations}
where the automorphism $\mathcal{T}_{i>k;j}$ occurs in the direction along which the quantity $\text{Im}[(\chi_k-\chi_j)/(\chi_j-\chi_i)]$ changes from negative to positive. The higher-order Stokes lines similarly define a graph $\mathcal{G}_H$ defined as the union $\bigcup_{i,j,k} h_{i>k;j}$. We note that higher-order Stokes lines pass through Stokes crossing points $\chi_i = \chi_j = \chi_k$ and emanate to and from turning points $\chi_i = \chi_j$ and $\chi_j = \chi_k$. The complete Stokes graph consisting of higher-order and ordinary Stokes lines is denoted by $\mathcal{G} = \mathcal{G}_S \cup \mathcal{G}_H$.

The ordinary Stokes phenomenon consists of a set of automorphisms $\mathcal{S}_{i>j} \subset GL_{N}(\mathbb{R})$ (when considered as acting on the vector of transseries parameters) associated to Stokes lines $l_{i>j}$. The higher-order Stokes automorphisms described above are now automorphisms of the Stokes automorphisms $\mathcal{S}_{i>j}$ and thus define elements of $\text{Aut}(GL_N(\mathbb{R}))$. The result is a representation of open curves $\gamma: [0,1] \to \mathbb{C}_z$ into the automorphism group acting on the transseries parameters $(\sigma_1,\ldots,\sigma_N)^T$. A curve may be decomposed into words of generators $\mathcal{S}_{i>j}$ in $\text{GL}_N(\mathbb{R})$ corresponding to Stokes lines, and generators $\mathcal{T}_{i>k;j}$ in $\text{Aut}(\text{GL}_N(\mathbb{R}))$ corresponding to higher-order Stokes automorphisms. In this way, the total Stokes phenomenon from $\gamma(0)$ to $\gamma(1)$ may be determined. Geometrically, the graph $\mathcal{G}$ induces a $GL_N(\mathbb{R}) \rtimes \text{Aut}(GL_N(\mathbb{R}))$ bundle over $\mathbb{C}$ with an associated parallel transport of the representation on the transseries parameters.

The activity of both regular Stokes lines $l_{i>j}$ and higher-order Stokes lines $h_{i>k;j}$ is defined in our framework by
\begin{equation}\label{eq:active}
l_{i>j} ~\text{is}~
\left\{\begin{aligned}
\text{active if}&\quad S_{ij} \neq 0,\\
\text{inactive if}&\quad S_{ij}=0,\\
\text{relevant if}&\quad  S_{ij}\sigma_i \neq 0,\\
\text{irrelevant if}&\quad  \sigma_i=0,
\end{aligned}\right.
\qquad
h_{i>k;j} ~\text{is}~
\left\{\begin{aligned}
\text{active if}\quad S_{ij}S_{jk} \neq 0,\\
\text{inactive if}\quad S_{ij}S_{jk}=0.
\end{aligned}\right.
\end{equation}
An intuitive description of the switching of activity of higher-order Stokes lines themselves may be given as follows.
Note that the higher-order Stokes automorphism $\mathcal{T}_{i>k;j}$ (from which the value of $S_{ik}$ changes) involves both $S_{ij}$ and $S_{jk}$. These have the potential to change value along $h_{i>k;j}$, which occurs when either $h_{i>j;l}$ or $h_{j>k;l}$ intersect with $h_{i>k;j}$. 
At least four different singulant functions are required for this to occur.
Crucially, a higher-order Stokes line will itself be inactive if either $S_{ij}$ or $S_{jk}$ is zero. 
Inactivity of higher-order Stokes lines can be interpreted from the perspective of the Borel plane (c.f. \cite{howls2004higher,shelton2025HOSP}), and also from adjacency graphs (which encode whether $S_{ij}$ is non-zero) arising in the steepest descent analysis of integrals \citep{berry1991hyperasymptotics}.
From the perspective of the Borel plane, inactivity due to $S_{ij}=0$ means that the Borel singularity $\chi_j-\chi_i$ is on a different Riemann sheet in the Borel plane.

\subsubsection{Stokes phenomenon at additional exponentials}
\label{sec:fiveormore}
The observation that a new Stokes phenomenon (the higher-order Stokes phenomenon) appears due to the interaction of three singulants leads to the natural question of whether new Stokes phenomenon will continue to appear when four or more singulants are present in the asymptotic transseries \eqref{eq:transseries0}.
The representation for $\psi_p^{(j)}$ in \eqref{eq:latetermsp} involved only the leading-order divergent term as $p\to \infty$, which had the amplitude function $\psi^{(k)}_{0}$.
By considering subdominant components of this with amplitude functions given by $\psi_{l}^{(k)}$ for instance, and then considering the late-late-late-terms of the expansion with $l \to \infty$, the result obtained is then identical to \eqref{eq:latetermsp}
but with $j \mapsto k$ and $p \mapsto l$.
Therefore, the same phenomenon is obtained.
Continuously reexpanding the amplitude functions simply produces the late-terms of a different component of the asymptotic transseries.

When four or more singulant functions are present in the transseries, the result of this reexpansions is that the higher-order Stokes automorphism $\mathcal{T}_{i>k;j}$ can change value along the pre-existing higher-order Stokes line $h_{i>k;j}$. Furthermore, this occurs at a higher-order Stokes crossing point where $h_{i>k;j}$ intersects with $h_{i>j;l}$ or $h_{j>k;l}$. Such phenomenon has previously been detected in the fifth-order PDE studied by \cite{body2005exponential}.
Therefore, no new Stokes phenomenon emerges occurs when transitioning from three to four WKBJ components in the asymptotic transseries, except for the fact that the higher-order automorphism can change value due to the crossing of higher-order Stokes lines.

The same observation is true when transitioning to a problem with five or more exponentials in the transseries, for which the only new effect is a more frequent appearance of crossing higher-order Stokes lines. Despite the increased number of exponential contributions, all changes in the higher-order Stokes automorphism occur across higher-order Stokes lines that describe an interaction between three exponential contributions.

\section{Asymptotic analysis of the swallowtail catastrophe}
\label{sec:swallowtail}
We now consider an example that contains each of the important features remarked upon in \S\,\ref{sec:method}, and which requires a careful study of both Stokes and higher-order Stokes automorphisms to fully determine the asymptotic behaviour.
This is the swallowtail integral, which is a generalisation of both the Airy and Pearcey integrals.
Stokes phenomenon in the Airy integral was famously studied in the pioneering work by \cite{stokes1864discontinuity},
and the Pearcey integral has subsequently been used to study the higher-order Stokes phenomenon associated with intersecting Stokes lines \citep{howls_2004,shelton2025HOSP}.
We demonstrate that the swallowtail integral contains crossing higher-order Stokes lines, which result in certain higher-order Stokes lines being inactive.
This phenomenon is resolved by applying our automorphisms for the Stokes constants of the asymptotic transseries, and a correct application of these higher-order Stokes automorphisms is necessary to enforce the quantization conditions on the problem.

In \S\,\ref{sec:swallowint} we convert the swallowtail integral into the singular perturbative form, derive the formal asymptotic transseries in \S\,\ref{sec:swallowtailtrans}, solve for the singulant functions in \S\,\ref{sec:singulantswallow}, and derive boundary conditions for the Stokes constants in \S\,\ref{sec:integralstokesconstant}.
The Stokes and higher-order Stokes automorphisms are then considered in \S\,\ref{sec:stokesstructure}.

\subsection{The swallowtail integral}
\label{sec:swallowint}
The swallowtail integral is a fundamental integral that arises in catastrophe theory, typically given in integral form as 
\begin{equation}
\label{eq:swallowint}
\Psi(x_1,x_2,x_3) = \int_{-\infty}^{\infty} \exp{\big(\mathrm{i}[t^5+x_3 t^3+x_2 t^2+x_1 t]\big)}\mathrm{d}t,
\end{equation}
This integral contains four distinct exponential contributions; it will serve as an example to demonstrate how the Stokes automorphism approach can track all of the possible Stokes phenomena that occur in such a system. Instead of the integral formulation \eqref{eq:swallowint}, we will focus on the swallowtail differential equation
\begin{equation}
\label{eq:swallowdiff}
5\frac{\mathrm{d}^4 \Psi}{\mathrm{d} x_1^4} -3 x_3\frac{\mathrm{d}^2 \Psi}{\mathrm{d} x_1^2} -2\mathrm{i}x_2\frac{\mathrm{d} \Psi}{\mathrm{d} x_1} +x_1 \Psi=0,
\end{equation}
where $x_1$, $x_2$, and $x_3$ are real-valued parameters, which is satisfied by the integral expression \eqref{eq:swallowint} for appropriate behavioural conditions. The Stokes line structure associated with solutions of \eqref{eq:swallowdiff} is independent of the specific choice of behavioural conditions which will be specified later in \eqref{eq:far-field}. 

We will develop asymptotic approximations of $\Psi$ in the limit of large $x_1$ by studying the Stokes phenomenon in the complex plane, $x_1 \in \mathbb{C}$.
It is possible to study the solution for $x_1 \in \mathbb{C}$ by analytically continuing both the integral \eqref{eq:swallowint} and the differential equation \eqref{eq:swallowdiff}, which were posed for $x_1 \in \mathbb{R}$.
We rescale the problem by introducing the small parameter $\epsilon$ through the definitions
\begin{equation}\label{eq:rescalings}
z= \epsilon^{4/5} x_1, \quad b=  \epsilon^{3/5} x_2, \quad a= \epsilon^{2/5} x_3, \quad \psi(z,a,b,\epsilon)=\Psi\bigg(\frac{z}{\epsilon^{4/5}},\frac{b}{\epsilon^{3/5}},\frac{a}{\epsilon^{2/5}}\bigg),
\end{equation}
for which the asymptotics of $\Psi$ for large $x_1$ correspond to that of $\psi$ for small $\epsilon$.
By substituting the rescalings \eqref{eq:rescalings} into the integral definition \eqref{eq:swallowint} and changing integration variables with $t \mapsto -\mathrm{i}t$ we obtain
\begin{equation}
\label{eq:newswallowint}
\psi(z,a,b,\epsilon) = \frac{1}{\mathrm{i} \epsilon^{1/5}}\int_{\mathscr{L}} \exp{\bigg(\frac{t^5-a t^3-\mathrm{i}b t^2+z t}{\epsilon}\bigg)}\mathrm{d}t.
\end{equation}
Under this change of variables, the contour $\mathscr{L}$ in integral \eqref{eq:newswallowint}, which takes values in $t \in \mathbb{C}$, starts at infinity in the sector $-3 \pi/10 <\text{arg}[t]<-\pi/10$ and ends at infinity in the sector $\pi/10 <\text{arg}[t]<3\pi/10$. 
The integral expression \eqref{eq:newswallowint} is a solution of the differential equation
\begin{subequations}
\label{eq:swallowtailmain}
\begin{equation}\label{eq:swallowtailrescaled}
    5\ep^4\frac{\de^4 \psi}{\de z^4} -3 a \ep^2\frac{\de ^2 \psi}{\de z^2} -2 \ii b\ep \frac{ \de \psi}{\de z}  +z \psi =0,
\end{equation}
which is obtained by substituting rescalings \eqref{eq:rescalings} into \eqref{eq:swallowdiff}, and has far-field asymptotics given by
\begin{equation}\label{eq:far-field}
\psi(z) \sim -\mathrm{i}\frac{\sqrt{\pi}}{\sqrt{2 \cdot 5^{1/4}\alpha^3 z^{3/4}}}\mathrm{e}^{-\frac{4 \alpha z^{5/4}}{5^{5/4} \epsilon}} \quad \text{as $|z|\to\infty$} ~\arg[z]=-3\pi/4,
\end{equation}
\end{subequations}
for $\alpha=\mathrm{e}^{ 3 \mathrm{i}\pi/4}$.
The far-field condition \eqref{eq:far-field} then ensures that the solution of ODE \eqref{eq:swallowtailrescaled} is equivalent to that given by the integral \eqref{eq:newswallowint}.
Since the integral definition will be found, via a steepest descent analysis in \S\,\ref{sec:integralstokesconstant}, to have transseries parameters $\sigma_1=1$ and $\sigma_2=\sigma_3=\sigma_4=0$ as $|z|\to \infty$ with $\arg[z]=-3\pi/4$, the far-field condition \eqref{eq:far-field} specifies that only one of the four possible WKBJ components of the asymptotic transseries is present in this region. The other three WKBJ components, which are not present in \eqref{eq:far-field}, have $\alpha=\mathrm{e}^{3 \mathrm{i}\pi/4}$ and $\alpha= \mathrm{e}^{\pm \mathrm{i}\pi/4}$. Note further that we have chosen a value of $\arg[z]$ in \eqref{eq:far-field} such that the specified exponential is the smallest out of all four permitted exponentials --- therefore since \eqref{eq:far-field} specifies it to be the dominant component of $\psi(z)$, the other three WKBJ components are unambiguously zero.

When solving for $\psi$ in the limit of $\epsilon \to 0$, we regard $z$ as the independent variable and both $a$ and $b$ as free parameters, the value of which will influence the resultant Stokes line structure.
Here, we focus on the case when $a$ and $b$ are real valued and $z$ is complex valued, but the method can easily be extended to the case where $a$ and $b$ take complex values.

\subsection{The asymptotic transseries}
\label{sec:swallowtailtrans}
We now solve for the asymptotic behaviour of $\psi(z,a,b,\epsilon)$ under the limit of $\epsilon \to 0$. For brevity, the functional dependence on $a$ and $b$ will be omitted from the notation.
The first step is to solve for each component of the asymptotic transseries \eqref{eq:transseries0}. Equations for the singulant and amplitude functions are obtained by substituting the transseries ansatz
\begin{equation}
\label{eq:swallowtailexp}
    \psi(z;\ep) \sim \sum_{i=1}^{4}\psi^{(i)}(z;\epsilon)\mathrm{e}^{-\chi_i(z)/\epsilon}
\end{equation}
into the swallowtail equation \eqref{eq:swallowtailrescaled}. By balancing powers of $\epsilon$, we find at leading order the equation
\begin{align}\label{eq:quartic}
(\chi_i^{\prime})^4-3a(\chi_i^{\prime})^2+2\ii b \chi_i^{\prime}+z=0,
\end{align}
and at the next order, 
\begin{equation}
\begin{aligned}\label{eq:ampequation}
&\ep^3 \frac{\mathrm{d}^4 \psi^{(i)}}{\mathrm{d}z^4}-20 \ep^2 \chi_i^{\prime}\frac{\mathrm{d}^3 \psi^{(i)}}{\mathrm{d}z^3}+3\Big[10 \ep (\chi_i^{\prime})^2-a \ep-10 \ep^2 \chi_i^{\prime \prime}\Big] \frac{\mathrm{d}^2 \psi^{(i)}}{\mathrm{d}z^2}+\Big[6 a \chi_i^{\prime}-2 \ii b -20 (\chi_i^{\prime})^3 \\
&+60 \ep \chi_i^{\prime}\chi_i^{\prime \prime}-20 \ep^2 \chi_i^{\prime \prime \prime}\Big]\frac{\mathrm{d} \psi^{(i)}}{\mathrm{d}z}+\Big[3a \chi_i^{\prime \prime}-30 \chi_i^{\prime \prime}(\chi_i^{\prime})^2+20 \ep \chi_i^{\prime}\chi_i^{\prime \prime \prime}+15 \ep (\chi_i^{\prime \prime})^2-5\ep^2 \chi_i^{\prime \prime \prime \prime}\Big]\psi^{(i)}=0.
\end{aligned}
\end{equation}
The leading order equation \eqref{eq:quartic} is a quartic equation for $\chi_i^{\prime}$, which may be solved analytically to find the four different exponents appearing in transseries \eqref{eq:swallowtailexp}.
Next, we expand the amplitude function $\psi^{(i)}(z;\epsilon)$ as $\epsilon \to 0$ by taking
\begin{equation}\label{eq:ampseries}
\psi^{(i)}(z;\ep) \sim \sum_{n=0}^{\infty}\ep^n \psi^{(i)}_n(z).
\end{equation}
Substituting \eqref{eq:ampseries} into \eqref{eq:ampequation} and balancing powers of $\epsilon$ yields at leading order as $\epsilon \to 0$
\begin{equation}
\label{eq:leadingorderamp}
\big[-20(\chi_i^{\prime})^3+6a \chi_i^{\prime}-2\ii b\big]\frac{\mathrm{d} \psi^{(i)}_0}{\mathrm{d}z} +\big[3a \chi _i^{\prime \prime}-30 \chi_i^{\prime \prime}(\chi_i^{\prime})^2 \big]\psi^{(i)}_0=0,
\end{equation}
which has the solution
\begin{equation}\label{eq:leadingamp}
\psi^{(i)}_0(z)=\frac{\sqrt{\pi}}{\sqrt{-10(\chi_i^{\prime})^3+3 a \chi_i^{\prime}-\ii b}}.
\end{equation}
The constant of integration has been set to be $\sqrt{\pi}$ to match with typical results obtained in the exact asymptotics of integrals, in which the Stokes constants take values $-1$, $0$, or $1$. Changing the value of this constant of integration subsequently changes the value of the transseries parameters $\sigma_i$ and Stokes constants $S_{ij}$, and so we are free to set these constants without any loss of generality. 

\subsubsection{Singular points}

Stokes lines are generated at points where the expansion \eqref{eq:swallowtailexp} fails to be asymptotic, which correspond to points where the amplitude \eqref{eq:leadingamp} becomes infinite. The leading-order amplitude function $\psi^{(i)}_0(z)$ is singular wherever the denominator in \eqref{eq:leadingamp} equals zero, which is any point at which the value of $\chi_i'$ satisfies $-10 (\chi_i')^3+3 a (\chi_i')-\mathrm{i}b=0$. 

By first solving the cubic equation $-10 (\chi_i')^3+3 a (\chi_i')-\mathrm{i}b=0$ for $t_i=\chi_i'$ at the singularities, we can find the values of $z$ at which the amplitude is singular by substituting $t_i$ into equation \eqref{eq:quartic} and solving directly for the associated values of $z$. This gives
\begin{equation}
\label{eq:singularity}
z_i=\frac{3 t_i}{2}\big( a t_i- \mathrm{i}b\big)
\end{equation}
for $i=1,2,3$. The asymptotic solution is therefore singular at the three locations $z_i \in \mathbb{C}$.

\begin{figure}
    \centering
    \includegraphics[scale=1]{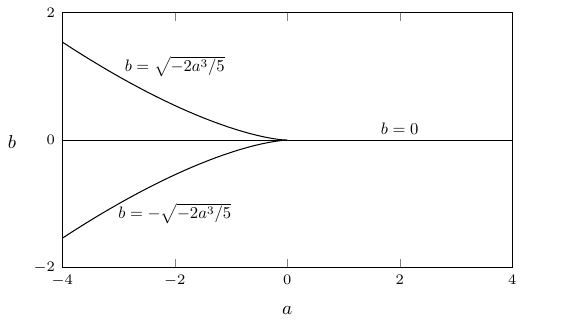}
    \caption{The values of $a$ and $b$ are shown that result in coalescence of the three singular points \eqref{eq:singularity}. Along $b=0$ two singularities coincide at $z=9a^2/20$, and on $a=-(5/2)^{1/3} |b|^{2/3}$ two singularities coincide at $z=-3a^2/20$. At $a=b=0$ all three singularities coincide at $z=0$. In each of these cases, the transseries \eqref{eq:swallowtailexp} still contains four distinct WKBJ components, but the Stokes geometry is simplified. The case of one singularity with $a=b=0$ contains no intersecting Stokes lines and so is analogous to the Stokes geometry for the Airy integral, and that with two singularities contains no crossing higher-order Stokes lines and so is analogous to the Stokes geometry for the Pearcey integral. The Stokes line structures for these simplified cases are shown later in figure~\ref{fig:otherab}.}
    \label{fig:abvalues}
\end{figure}

In general, $z_1$, $z_2$, and $z_3$ are distinct points.
However, for certain values of $a$ and $b$, the value of these singular points coincide with one another along the real axis.
For instance, when $a=b=0$ all three singular points are at $z=0$.
When $b=0$ there are two distinct singular points for all values of $a$, one at $z=0$ and another repeated singularity at $z=9a^2 /20$.
Additionally, only two singular points emerge for $b \neq 0$ when $a=-(5/2)^{1/3}|b|^{2/3}$, and these are given by $z=6a^2/5$ and $z=-3a^2/20$, with the latter being the repeated singular point.
These values of $a$ and $b$ are shown in figure~\ref{fig:abvalues}.
Note that all of these cases result in a transseries expansion with four WKBJ components, but that with two singularities contains no crossing higher-order Stokes lines, and that with one singularity contains no crossing Stokes lines. 
In \S\,\ref{sec:stokesstructure}, results are given for the values $(a,b)=(1,3)$ which produces an expansion with three distinct singular points for $z \in \mathbb{C}$, and a Stokes line structure that is symmetric about $\text{Im}[z]=0$. 
This example will contain the full range of Stokes and higher-order Stokes automorphisms that can emerge, in which a Stokes constant changes along a higher-order Stokes line due to intersection with another higher-order Stokes line.
In the discussion of \S\,\ref{sec:conc} we briefly illustrate the Stokes structure for the simpler cases of $(a,b)=(0,0)$ and $(a,b)=(1,0)$.

\subsection{The singulant functions}
\label{sec:singulantswallow}
Four different solutions for $\chi_i^{\prime}(z)$ are obtained from the quartic equation \eqref{eq:quartic}.
These must be integrated in order to analyse the resultant Stokes line structure.
We now integrate these analytically by solving the Hamilton-Jacobi equation that arises from the Borel operator of the differential equation \citep{kawai2005algebraic}.
The Borel operator, $\mathscr{P}_B$, and symbol, $\sigma$, of this equation are given by
\begin{equation}\label{eq:Appborelop}
\mathscr{P}_B= 5\partial_z^4 -3a \partial^2_w \partial^2_z - 2\ii b \partial^3_w \partial_z + z \partial^4_w \qquad \text{and} \qquad \sigma = 5\xi^4 -3a \eta^2 \xi^2 - 2\ii b \eta^3 \xi + z \eta^4,
\end{equation}
where $\xi$ and $\eta$ are local coordinates on the cotangent directions.
The Hamiltonian equations and corresponding boundary conditions arising from the complex symplectic structure on $T^*(\mathbb{C}_w \times \mathbb{C}_z)$ are given by
\begin{subequations}
\begin{equation}
\label{eq:hamiltonianequations}
 \dot{z}(\tau) = \partial_\xi \sigma, 
 \quad \dot{w}(\tau) = \partial_{\eta} \sigma, 
 \quad \dot{\xi}(\tau) = -\partial_{z}\sigma, 
 \quad  \dot{\eta}(\tau) = -\partial_w \sigma, 
 \quad  \sigma (\tau) = 0,
\end{equation}
\begin{equation}
 z(0)=0, \quad w(0)=w_0, \quad \xi(0)=0, \quad \eta(0)=1, \quad \sigma(0)=0.
\end{equation}
\end{subequations}
The equations for $\eta$ and $\xi$ may be solved directly to find $\eta(\tau)=1$ and $\xi(\tau)=-\tau$, which we use to solve the remaining equations for $z$ and $w$, yielding
\begin{subequations}
\refstepcounter{equation}\label{eq:zsol}
\refstepcounter{equation}\label{eq:wsol}
\begin{equation}
{z}(\tau) =-5 \tau^4 +3a \tau^2 - 2 \mathrm{i}b \tau, \qquad 
{w}(\tau) =-4 \tau^5 + 2 a \tau^3 -\mathrm{i}b \tau^2 +w_0.
\tag{\ref*{eq:zsol},$b$}
\end{equation}
\end{subequations}
Equation \eqref{eq:zsol} is recognised as the same quartic equation that governs $\chi_i^{\prime}(z)$ in \eqref{eq:quartic}, and so $\tau=\chi_i^{\prime}(z)$.
Furthermore, by evaluation of the chain rule on $\mathrm{d} w/ \mathrm{d}z = (\mathrm{d}w/\mathrm{d}{\tau}) (\mathrm{d} {\tau} / \mathrm{d} z)$, it can be shown that $\mathrm{d}w/\mathrm{d}z=\mathrm{d}\chi / \mathrm{d}z$.
This implies that $w=\chi(z)$ (up to a constant of integration), and so from \eqref{eq:wsol} we have found an exact solution for each singulant, $\chi_i(z)$, in terms of its derivative, $\chi_i^{\prime}(z)$.
In substituting for the definition of $(\chi_i^{\prime})^4$ from the quartic \eqref{eq:quartic}, this exact solution is given by
\begin{subequations}
\begin{equation}\label{eq:singsol}
\chi_i(z)=-\frac{2a}{5}(\chi_i^{\prime})^3 +\frac{3 \mathrm{i}b}{5}(\chi_i^{\prime})^2+\frac{4z}{5}\chi_i^{\prime}+w_0
\end{equation}
for $i=1,2,3,4$, which satisfies the boundary conditions
\begin{equation}\label{eq:singBC}
\chi_2(z_1)-\chi_1(z_1)=0, \qquad  \chi_4(z_2)-\chi_3(z_2)=0, \qquad \chi_3(z_3)-\chi_2(z_3)=0,
\end{equation}
\end{subequations}
which hold at each of the three turning points.
Thus, the constant $w_0$ is arbitrary and is determined by the far-field conditions of the problem which were given in \eqref{eq:far-field} and yield $w_0=0$.

Note that it is also possible to derive the solution \eqref{eq:singsol} from the original integral formulation \eqref{eq:newswallowint}.
Here, the exponent of the saddle point integral is $f(t,z) = t^5 - at^3 - ib t^2 + zt$, and so the saddle point locations in the complex $t$-plane are solutions of $\partial_tf(t_i,z) =0$.
There are four saddle point locations $t_i$ in total, which are given as solutions to $5t_i^4 - 3at_i^2 - 2ib t_i + z = 0$ for $i=1,2,3,4$.
The singulants are then defined by $\chi_i(z) = -f(t_i,z)$.
On the other hand we note that $ \partial_z f(z,t) = t$, so that $t_i=-\chi_i'(z)$.
Then we immediately observe the relation $\chi_i(z) = [2a(t_i)^3 +3 \mathrm{i}b(t_i)^2-4z t_i]/5$,
which matches with \eqref{eq:singsol}.
While obtaining a solution for $\chi_i(z)$ from the integral formulation is the easier of the two methods given in this section, solving the Hamiltonian system is the more general approach which can be applied to linear differential equations.

\subsection{Determination of the Stokes constants}
\label{sec:integralstokesconstant}

The Stokes constants $S_{ij}$ are analytic invariants of the underlying differential equation \eqref{eq:swallowtailrescaled}.
The value of these will change across $\mathbb{C}_z$ on account of the higher-order Stokes automorphisms \eqref{eq:HOSA+P}.
We now determine the value of these at the origin, $z=0$.

Multiple methods exist to calculate the values of the Stokes constants $S_{ij}$.
For instance, the formula \eqref{eq:Stokesconstantform} yields values for $S_{ij}$ in different sectors of $\mathbb{C}_z$, and the method of hyperterminants \citep{bennett2018globally} produces values for every Stokes constant at the same location in $\mathbb{C}_z$.
Here, we will use the integral solution \eqref{eq:newswallowint} to solve for the Stokes constants by the method of \cite{berry1991hyperasymptotics}.
We first consider the integral representations for both the expansion $\psi^{(i)}(z;\epsilon) \sim \psi_0^{(i)}(z)+ \epsilon \psi_1^{(i)}(z)+\cdots$ and $\psi_n^{(i)}$ for $n \geq 0$, which are given by
\begin{subequations}
\begin{align}
\label{eq:intseriessol}
\sum_{n=0}^{\infty} \epsilon^n \psi_n^{(i)}(z) & \sim \frac{1}{\epsilon^{1/2}}\int_{C_i(\text{arg}[\epsilon])} \exp\bigg( -\frac{f(t,z)-\chi_i(z)}{\epsilon} \bigg) \mathrm{d} t,\\
\label{eq:intlateterm}
\psi_n^{(i)}(z)&=\frac{1}{2 \pi \mathrm{i}} \ointctrclockwise \frac{\Gamma(n+1/2)}{[f(t,z)-\chi_i(z)]^{n+1/2}} \mathrm{d}t.
\end{align}
\end{subequations}
The path of integration in \eqref{eq:intseriessol} is the steepest descent contour $C_i(\text{arg}[\epsilon])$ which passes through the saddle point $t=\tilde{t}$ at which $\partial_{t} f (\tilde{t},z)=0$, and that in \eqref{eq:intlateterm} is the positively oriented contour enclosing the path of steepest descent through $\tilde{t}$.

Recall that the Stokes constants $S_{ij}$ appear in the late-term expansion \eqref{eq:latetermsall} of $\psi_n^{(i)}(z)$ as $n \to \infty$.
Therefore, we need to manipulate \eqref{eq:intlateterm} for large $n$ by deforming the positively oriented contour of integration onto other paths of steepest descent. 
Specifically, we consider the values of $\text{arg}[\epsilon]$ that satisfy $\text{arg}[(\chi_j-\chi_i)/\epsilon]=0$, and denote these values by $\text{arg}[\epsilon]=-\sigma_{ij}$ (where $j \neq i$). 
Then, for a specified value of $i$ there are at most three steepest descent contours $C_j(\sigma_{ij})$ that the positively oriented contour from \eqref{eq:intlateterm} can be deformed onto, and this yields the result
\begin{equation}\label{eq:intlateterm4}
\psi_n^{(i)}(z)\sim\frac{\mathrm{i}}{2 \pi} \sum_{j, j\neq i} (-1)^{\gamma_{ij}}A_{ij} \sum_{p=0}^{\infty} \psi_p^{(j)}(z)\frac{\Gamma(n-p)}{[\chi_j(z)-\chi_i(z)]^{n-p}}.
\end{equation}
Here, $A_{ij}$ is the adjacency constant which takes a value of zero if the deformation of the contour from \eqref{eq:intlateterm} does not contain the contour $C_{j}(\sigma_{ij})$, and one otherwise. The constant $\gamma_{ij}$ is equal to zero if the deformed contour is aligned with $C_j(\sigma_{ij})$, and is equal to one if the deformed contour is in an opposite direction to $C_j(\sigma_{ij})$.
By comparing \eqref{eq:intlateterm4} for the late-terms of the Swallowtail problem to the general late-term expansion \eqref{eq:latetermsall},
we obtain the solution
\begin{equation}
\label{eq:stokesconstantsgamma}
   S_{ij}=(-1)^{\gamma_{ij}} A_{ij}
\end{equation}
for the Stokes constants. 
Since the adjacency constant $A_{ij}$ is either $0$ or $1$, we have that the Stokes constants $S_{ij}$ take values of either $-1$, $0$, or $1$.
Note also that in this specific example, due to the Swallowtail equation having a simple integral solution and also due to the choice of integration in the leading-order amplitude function \eqref{eq:leadingamp}, we have that $S_{ij}=-S_{ji}$.

\begin{figure}
    \centering
    \includegraphics[scale=1]{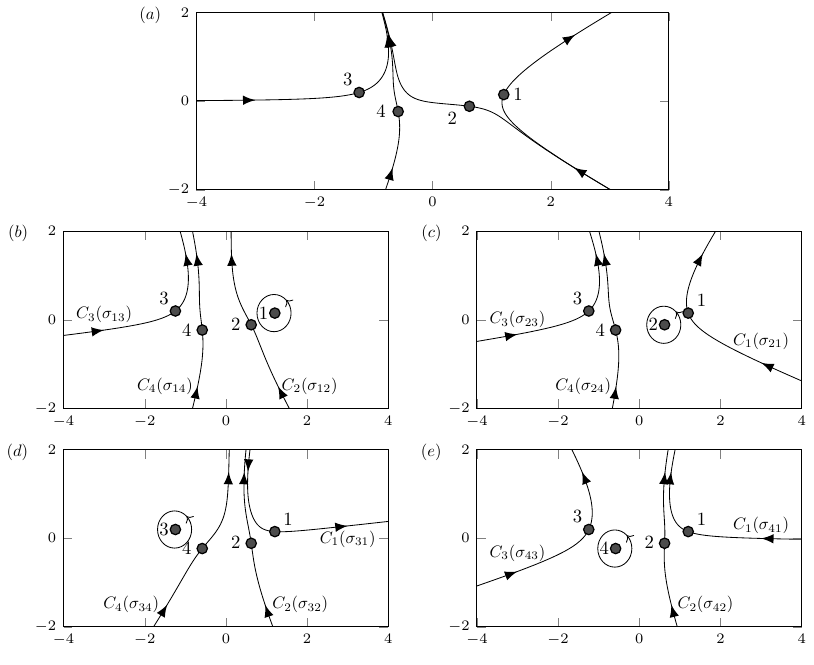}
    \caption{The paths of steepest descent $C_i(\arg[\epsilon])$, defined as where $(f(t,z)-\chi_i)/\epsilon$ is real and positive emerging from saddle $i$, are shown for the swallowtail integral \eqref{eq:newswallowint} at $z=3+0.5\mathrm{i}$. $(a)$ shows the steepest descent curves $C_i(0)$ through each saddle for the case where $\arg[\epsilon]=0$. The arrows along each of these define the direction of integration in \eqref{eq:intseriessol}.
    In $(b1$--$b2)$, the special steepest descent curves $C_j(\sigma_{ij})$ are shown for $\sigma_{ij}=\arg[\chi_j-\chi_i]$. The direction along each steepest descent curve is obtained by taking the direction in $(a)$ and changing $\arg[\epsilon]$ continuously from $0$ to $\sigma_{ij}$. Values for the Stokes constants $S_{ij}$ are obtained from \eqref{eq:stokesconstantsgamma} by deforming the anticlockwise contour around saddle $i$ in panel $(bi)$ onto all adjacent curves $C_{j}(\sigma_{ij})$.}
    \label{fig:stokesconstants}
\end{figure}
Since the choice of index, from $i=1$ to $i=4$, in the singulant functions $\chi_i(z)$ is arbitrary, we specify this with the values at $z=3+0.5\mathrm{i}$ of $\chi_1=-1.0464 + 0.7948\mathrm{i}$, $\chi_2=-1.1944 + 0.0920\mathrm{i}$, $\chi_3=1.2212 + 2.0196\mathrm{i}$, and $\chi_4=1.0196 + 0.6936\mathrm{i}$, given to four decimal places. 
The steepest descent contours $C_{i}(\arg[\epsilon])$ are shown in figure~\ref{fig:stokesconstants}.
The arbitrary direction of integration along these in integral \eqref{eq:intseriessol} is specified by the arrows in figure~\ref{fig:stokesconstants}$(a)$.
Then, to work out the orientation of the special steepest descent curves $C_{j}(\sigma_{ij})$, as required for formula \eqref{eq:stokesconstantsgamma}, we start with the specified orientation for $\arg[\epsilon]=0$, and change this value continuously to $\sigma_{ij}$.
The resultant orientation of the special steepest descent curves are shown in panels $(b1)$--$(b4)$ of figure~\ref{fig:stokesconstants}.
By deforming the anticlockwise contour about each saddle point onto all adjacent steepest curves, we find that at $z=3+0.5 \mathrm{i}$ the Stokes constants $S_{ij}$ are given by
\begin{equation}
\label{eq:stokesfirstregion}
\begin{aligned}
S_{12}=-1, \quad S_{13}=0, \quad S_{14}=0, \qquad S_{21}=1, \quad S_{23}=0, \quad S_{24}=-1, \\
S_{31}=0, \quad S_{32}=0, \quad S_{34}=1, \qquad S_{41}=0, \quad S_{42}=1, \quad S_{43}=-1. 
\end{aligned}
\end{equation}
In the next section, we will compute the higher-order Stokes automorphisms that cause the values of $S_{ij}$ reported here for $z=3+0.5 \mathrm{i}$ to change across higher-order Stokes lines.

\section{The swallowtail Stokes line structure}
\label{sec:stokesstructure}
We now calculate the set of higher-order Stokes lines $h_{i>k;j}$ and the set of Stokes lines $l_{i>j}$ for the swallowtail system \eqref{eq:swallowtailmain}.
This is performed by using the singulant functions obtained in \eqref{eq:singsol}, as well as conditions 
\eqref{eq:HOSA+P} for $h_{i>k;j}$, and \eqref{eq:level1line} for $l_{i>j}$.
We will also need to calculate the Stokes constants $S_{ij}$ throughout the complex plane by starting with the values obtained in solution \eqref{eq:stokesfirstregion} at $z=3+0.5\mathrm{i}$ and then applying the higher-order Stokes automorphism $\mathcal{T}_{i>k;j}$ from \eqref{eq:HOSA+P2}.
Once the Stokes constants are known throughout $\mathbb{C}_{z}$, we then solve for the set of transseries parameters $\sigma_i$ by starting from values obtained from the far-field conditions \eqref{eq:far-field}.
\begin{figure}
    \centering
    \includegraphics[scale=1]{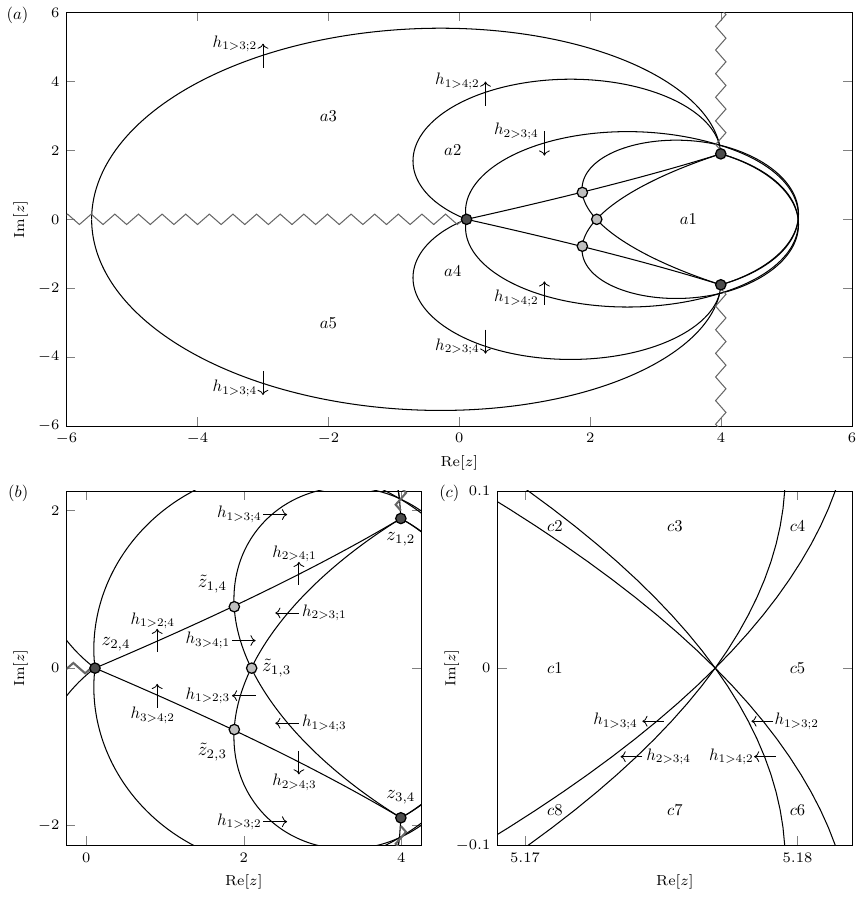}
    \caption{$(a)$ The higher-order Stokes line geometry for the swallowtail problem \eqref{eq:swallowtailrescaled} is shown for the parameter values $(a,b)=$ $(1,3)$. Turning points are shown with black circles, virtual turning points with grey circles, branch cuts with grey lines, and higher-order Stokes lines with black lines. Arrows across $h_{i>k;j}$ show the direction to which the automorphism $\mathcal{T}_{i>k;j}$ is applied. In $(b)$, the region near the turning points is shown. In $(c)$, the region with four intersecting higher-order Stokes lines is shown. Note that both $h_{i>k;j}$ and $h_{k>i;j}$ coincide along the same curve, but that the automorphisms are applied in opposite directions.}
    \label{fig:HOSL}
\end{figure}

\begin{figure}
    \centering
    \includegraphics[scale=1]{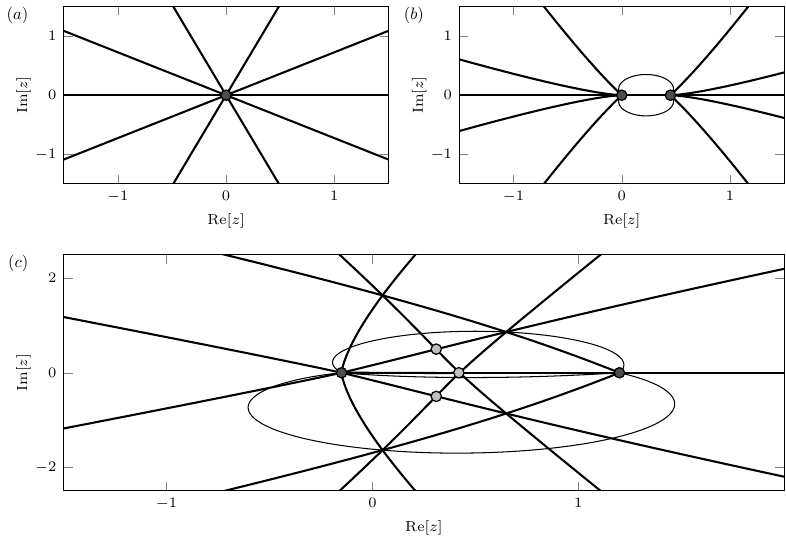}
    \caption{The Stokes lines $l_{i>j}$ (bold) and higher-order Stokes lines $h_{i>k;j}$ (thin) are shown for three cases of $(a,b)=(0,0)$ in $(a)$, $(a,b)=(1,0)$ in $(b)$, and $(a,b)=(-1,\sqrt{2/5})$ in $(c)$. These values of $a$ and $b$ result in either one or two turning points as shown in figure~\ref{fig:abvalues}. Stokes lines emanate from both the turning points (dark grey circles) and the virtual turning points (light grey circles).}
    \label{fig:otherab}
\end{figure}

Recall that in addition to $\epsilon$, the governing differential equation \eqref{eq:swallowtailrescaled} contains two parameters $a$ and $b$, for which we choose here the values $a=1$ and $b=3$. 
This yields three distinct singular points of the asymptotic expansion, which is the simplest setup that will result in a change in activity of the higher-order Stokes lines. 
Note that for certain values of $a$ and $b$, shown previously in figure~\ref{fig:abvalues}, fewer singular points emerge, and so the higher-order Stokes lines do not intersect.
In these special cases, shown later in figure~\ref{fig:otherab}, the asymptotic phenomenon that emerges is analogous to that found in the third-order Pearcey problem \citep{howls_2004,shelton2025HOSP}.

All higher-order Stokes lines $h_{i>k;j}$, both active and inactive, are shown in figure~\ref{fig:HOSL} for $(a,b)=(1,3)$.
The arrows across each higher-order Stokes line indicate the direction in which the higher-order Stokes automorphism $\mathcal{T}_{i>k;j}$ from \eqref{eq:HOSA+P2} is applied.
In addition to the three turning points $z_{i,j}$, at which $\chi_i=\chi_j$, there are also three virtual turning points $\tilde{z}_{i,j}$ at which the same condition holds. 
The difference between $z_{i,j}$ and $\tilde{z}_{i,j}$ is that the virtual turning points do not correspond to actual singularities that manifest in the asymptotic expansion due to the corresponding Stokes constant being zero nearby.
The higher-order Stokes lines emanating from virtual turning points $\tilde{z}_{i,j}$, as shown in figure~\ref{fig:HOSL}$(b)$, are therefore inactive.
Note that the branch cuts, shown with grey zig-zag lines, are branch cuts in each WKBJ component of the asymptotic solution, rather than a branch cut in the overall function. Such asymptotic branch cuts result in a reindexing of the Stokes and higher-order Stokes lines when crossed.
In figure~\ref{fig:HOSL}$(c)$, the region with four intersecting higher-order Stokes lines is shown in detail. Multiple Stokes lines will be seen to pass through the same crossing point and so values for the Stokes constants $S_{ij}$ are required in each sector here.

\begin{table}
\begin{tabular}{ cc }
    \begin{tabular}[t]{c||c|c|c|c|c|c}
        Region & $S_{12}$ & $S_{13}$ & $S_{14}$ & $S_{23}$ & $S_{24}$ & $S_{34}$\\
        \hline
        $a1$ &$-1$ & 0 & 0 & 0 & $-1$ & $1$ \\
        $a2$ & $-1$ & $0$ & $0$ & $-1$ & $-1$ & $1$ \\
        $a3$ & $-1$ & $0$ & $1$ & $-1$ & $-1$ & $1$ \\
        $a4$ & $-1$ & $0$ & $-1$ & $0$ & $-1$ & $1$ \\
        $a5$ & $-1$ & 0 & $-1$ & $1$ & $-1$ & $1$ \\
    \end{tabular} &
       \begin{tabular}[t]{c||c|c|c|c|c|c}
        Region & $S_{12}$ & $S_{13}$ & $S_{14}$ & $S_{23}$ & $S_{24}$ & $S_{34}$\\
        \hline
        $c1$ & $-1$ & 0 & 0 & 0 & $-1$ & $1$ \\
        $c2$ & $-1$ & 0 & 0 & 0 & $-1$ & $1$  \\
        $c3$ & $-1$ & 0 & $-1$ & 0 & $-1$ & $1$ \\
        $c4$ & $-1$ & 0 & $-1$ & $-1$ & $-1$ & $1$\\
        $c5$ &  $-1$ & $-1$ & $-1$ & $-1$ & $-1$ & $1$\\
        $c6$ &  $-1$ & $0$ & $-1$ & $-1$ & $-1$ & $1$ \\
        $c7$ & $-1$ & $0$ & $0$ & $-1$ & $-1$ & $1$ \\
        $c8$ & $-1$ & $0$ & $0$ & $0$ & $-1$ & $1$ \\
    \end{tabular}
    \end{tabular}
    \caption{Values for the Stokes constants $S_{ij}$ are shown for each of the regions from figure~\ref{fig:HOSL}. Note that $S_{ij}=-S_{ji}$ for the Swallowtail problem due to the choice of integration constant in the leading-order amplitude function \eqref{eq:leadingorderamp}. There are twelve Stokes constants in total, and due to this symmetry only six are reported here. Values for the Stokes constants in regions $a1$ and $c1$ were obtained in solution \eqref{eq:stokesfirstregion}, and the values in other regions obtained by applying the higher-order Stokes automorphism $\mathcal{T}_{i>k;j}$ from \eqref{eq:HOSA+P2}. The adjacency graph in each region of $\mathbb{C}_z$ is constructed by connecting nodes $i$ and $j$ by an edge if $S_{ij}$ contains a nonzero entry for that region.}
    \label{tab:Sij}
\end{table}
Values for the Stokes constants $S_{ij}$ across the complex plane are given in table~\ref{tab:Sij}. Here, the five regions $a1$--$a5$ are those labelled in figure~\ref{fig:HOSL}$(a)$, and the eight regions $c1$--$c8$ are those labelled in figure~\ref{fig:HOSL}$(c)$.
The values for $S_{ij}$ in region $a1$ and $c1$ were obtained in \S\,\ref{sec:integralstokesconstant} by studying the integral formulation of the problem, and the resultant values in other regions of $\mathbb{C}_z$ are obtained from these by applying the automorphism $\mathcal{T}_{i>k;j}$ from \eqref{eq:HOSA+P2} across each higher-order Stokes line shown in figure~\ref{fig:HOSL}.

\begin{figure}
    \centering
    \includegraphics[scale=1]{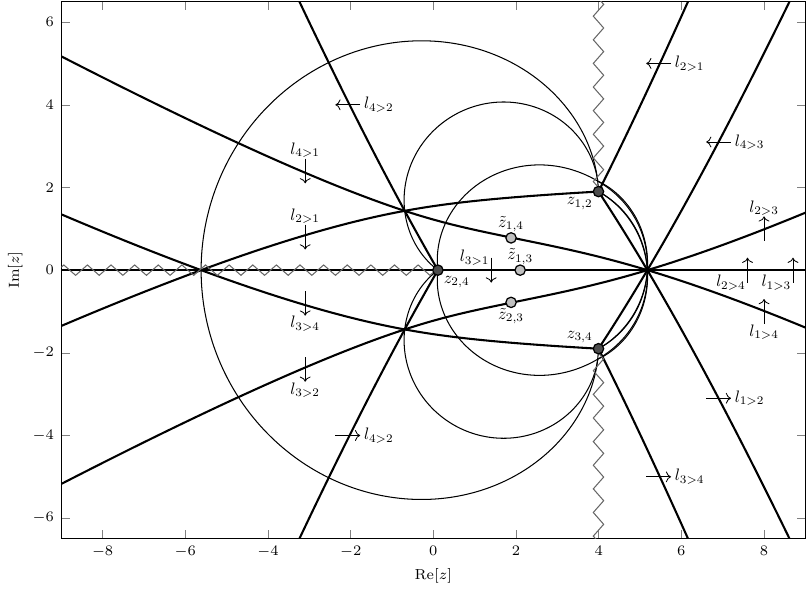}
    \caption{Stokes lines $l_{i>j}$ (bold lines) and higher-order Stokes lines $h_{i>k;j}$ (thin lines) are shown for the swallowtail problem \eqref{eq:swallowtailmain}. The Stokes lines $l_{i>j}$ satisfy criteria \eqref{eq:level1line} with the singulant functions \eqref{eq:singsol}. Note that while every higher-order Stokes line was shown figure~\ref{fig:HOSL}, only active higher-order Stokes lines are shown here. The higher-order Stokes line $h_{i>k;j}$ is inactive if the automorphism \eqref{eq:HOSA+P2} is suppressed due to either $S_{ij}=0$ or $S_{jk}=0$.}
    \label{fig:HOSLandSL}
\end{figure}
Since we have that $S_{13}=0$, $S_{14}=0$, and $S_{23}=0$ in region $a1$, all of the higher-order Stokes lines labelled in figure~\ref{fig:HOSL}$(b)$ that emanate from the virtual turning points are inactive.
The resultant active higher-order Stokes lines are shown in figure~\ref{fig:HOSLandSL} alongside every Stokes line $l_{i>j}$.
Here, Stokes lines are seen to intersect with one another in four locations.
The three of these locations with $\text{Re}[z]<0$ have three crossing Stokes lines, and the intersection point along the positive real axis has six intersecting Stokes lines.
The activity of certain Stokes lines will change at each of these Stokes crossing points due to a change in value of the Stokes constant $S_{ij}$.
For instance, the Stokes constant $S_{14}$ changes from a value of $0$ in region $a2$ to a value of $1$ in region $a3$. This results in a change of activity of the Stokes line $l_{4>1}$ either side of the respective Stokes crossing point.
Further, near the crossing point along the positive real axis the three Stokes constants $S_{13}$, $S_{14}$, and $S_{23}$ are all zero on the left (region $c1$ in table~\ref{tab:Sij}), and non-zero on the right of the crossing point (region $c5$ in table~\ref{tab:Sij}).
This results in the change of activity of three Stokes lines through this second Stokes crossing point.

\begin{figure}
    \centering
    \includegraphics[scale=1]{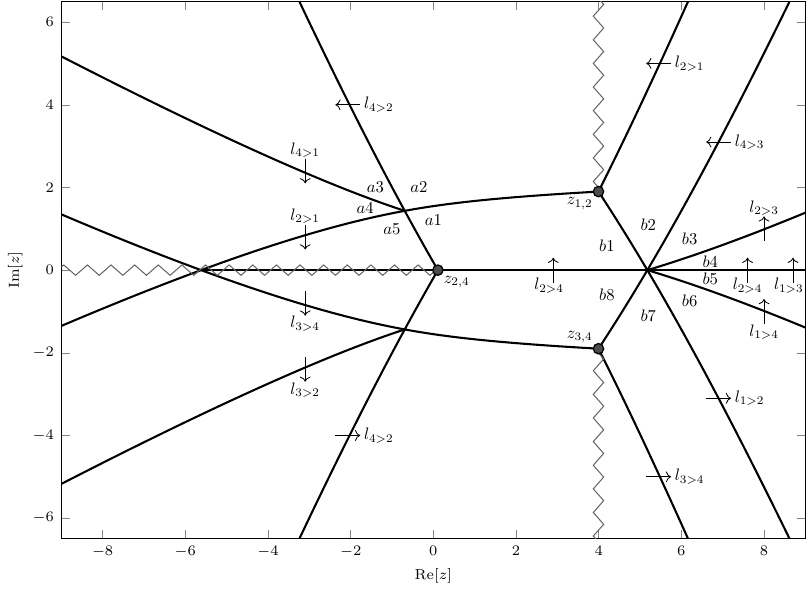}
    \caption{The active $l_{i>j}$ Stokes lines are shown for the swallowtail problem \eqref{eq:swallowtailmain}. These were obtained by taking all the Stokes lines $l_{i>j}$ shown in figure~\ref{fig:HOSLandSL} and retaining only those with $S_{ij} \neq 0$ along them. Values for the Stokes constants $S_{ij}$ are taken from table~\ref{tab:Sij}.}
    \label{fig:activeSL}
\end{figure}
The active Stokes lines $l_{i>j}$, which have $S_{ij}\neq 0$ along them, are shown in figure~\ref{fig:activeSL}.
We see that for the Stokes crossing point on the positive imaginary axis, the three Stokes lines $l_{2>3}$, $l_{1>3}$, and $l_{1>4}$ change activity thought the crossing point.
Values for the transseries parameters $\sigma_i$ in sectors surrounding each of the two labelled Stokes crossing points in figure~\ref{fig:activeSL} are given in table~\ref{tab:sigmai}.
To calculate these, we have computed the Stokes automorphisms \eqref{eq:level1auto} across each Stokes line $l_{i>j}$ for which the Stokes constants $S_{ij}$ from table~\ref{tab:Sij} have been used.
Note that while the Stokes constants $S_{ij}$ are fixed for a given differential equation, the values of the transseries parameters remain arbitrary from the equation and are instead set by boundary conditions.
We have therefore left these unspecified by starting with values $\sigma_i=\beta_i$ in the region labelled by both $a1$ and $b1$.

\begin{table}
\begin{tabular}{ cc }
    \begin{tabular}[t]{c||c|c|c|c}
         & $\sigma_{1}$ & $\sigma_{2}$ & $\sigma_{3}$ & $\sigma_{4}$\\
        \hline
        $a1$ &$\beta_1$ & $\beta_2$ & $\beta_3$ & $\beta_4$ \\
        $a2$ & $\beta_1-\beta_2$ & $\beta_2$ & $\beta_3$ & $\beta_4$ \\
        $a3$ & $\beta_1-\beta_2$ & $\beta_2+\beta_4$ & $\beta_3$ & $\beta_4$ \\
        $a4$ & $\beta_1-\beta_2-\beta_4$ & $\beta_2+\beta_4$ & $\beta_3$ & $\beta_4$ \\
        $a5$ & $\beta_1$ & $\beta_2+\beta_4$ & $\beta_3$ & $\beta_4$ \\
              $a1$ &$\beta_1$ & $\beta_2$ & $\beta_3$ & $\beta_4$ \\
    \end{tabular} &
       \begin{tabular}[t]{c||c|c|c|c}
        & $\sigma_{1}$ & $\sigma_{2}$ & $\sigma_{3}$ & $\sigma_{4}$\\
        \hline
        $b1$ &$\beta_1$ & $\beta_2$ & $\beta_3$ & $\beta_4$ \\
        $b2$ & $\beta_1$ & $\beta_2-\beta_1$ & $\beta_3$ & $\beta_4$  \\
        $b3$ & $\beta_1$ & $\beta_2-\beta_1$ & $\beta_3+\beta_4$ & $\beta_4$ \\
        $b4$ & $\beta_1$ & $\beta_2-\beta_1$ & $\beta_3+\beta_4+\beta_2-\beta_1$ & $\beta_4$\\
        $b5$ &  $\beta_1$ & $\beta_2-\beta_1$ & $\beta_3+\beta_4+\beta_2$ & $\beta_4+\beta_2-\beta_1$\\
        $b6$ &  $\beta_1$ & $\beta_2-\beta_1$ & $\beta_3+\beta_4+\beta_2$ & $\beta_4+\beta_2$ \\
        $b7$ & $\beta_1$ & $\beta_2$ & $\beta_3+\beta_4+\beta_2$ & $\beta_4+\beta_2$\\
        $b8$ & $\beta_1$ & $\beta_2$ & $\beta_3$ & $\beta_4+\beta_2$ \\
         $b1$ &$\beta_1$ & $\beta_2$ & $\beta_3$ & $\beta_4$ \\
    \end{tabular}
    \end{tabular}
    \caption{Values for the transseries parameters $\sigma_i$ are given in regions surrounding two Stokes crossing points shown in figure~\ref{fig:activeSL}. The first Stokes crossing point, which arises due to the intersection of three Stokes lines, divides the local domain into five sectors, labelled $a1$--$a5$. The second, which arises due to the intersection of six Stokes lines, divides the local domain into eight sectors, labelled $b1$--$b8$.}
    \label{tab:sigmai}
\end{table}

\begin{figure}
    \centering
    \includegraphics[scale=1]{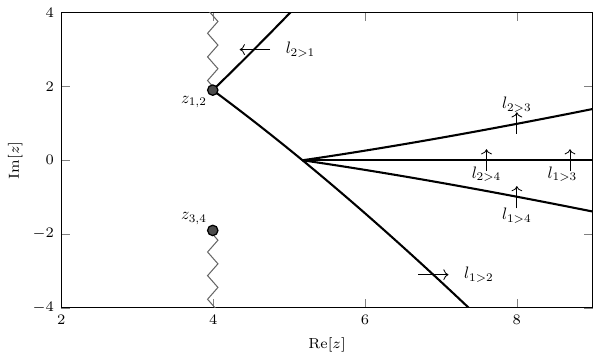}
    \caption{The relevant $l_{i>j}$ Stokes lines are shown for the swallowtail problem \eqref{eq:swallowtailmain}. These were obtained by taking all the active Stokes lines $l_{i>j}$ shown in figure~\ref{fig:activeSL} and retaining only those with $\sigma_i \neq 0$ along them. Values for the transseries parameters $\sigma_{i}$ are taken from table~\ref{tab:sigmai} with the specification of $\beta_1=1$, and $\beta_2=\beta_3=\beta_4=0$ which arises from the behavioural condition \eqref{eq:far-field}.}
    \label{fig:relevantSL}
\end{figure}

The specific choice arising from far-field condition \eqref{eq:far-field}, which was specified to match with the boundary conditions implicit in the swallowtail integral \eqref{eq:newswallowint}, corresponds to 
\begin{equation}
\label{eq:specificsigma}
\beta_1=1, \quad \beta_2=0, \quad \beta_3=0, \quad \beta_4=0.
\end{equation}
Therefore, many of the active Stokes lines shown in figure~\ref{fig:activeSL} are irrelevant since the corresponding exponential that induces Stokes phenomenon is not present.
Recall from \eqref{eq:active} that the Stokes line $l_{i>j}$ is relevant if both $\sigma_i\neq 0$ and $S_{ij}\neq 0$ along it.
The relevant Stokes lines for the Swallowtail problem with behavioural condition \eqref{eq:far-field}, which results in the transseries parameters from table~\ref{tab:sigmai} with values specified by \eqref{eq:specificsigma}, are shown in figure~\ref{fig:relevantSL}.

\section{Discussion and further directions}
\label{sec:conc}
In this study, we considered a complex function whose asymptotic behaviour is represented as a formal asymptotic transseries. The asymptotic behaviour can be completely characterised using two sets of fundamental constants.
These are the transseries parameters $\sigma_i$, which control the amplitude of each exponential, and the Stokes constants $S_{ij}$, which dictate the amount of the early terms of one exponential present in the late-terms of another. The Stokes constants do not depend on boundary data or far-field conditions, and characterise all solutions of the differential equation.
The transseries parameters are set by boundary data or far-field conditions of the problem, and change value throughout the complex plane on account of the Stokes phenomenon.

We interpreted the Stokes phenomenon as an automorphism \eqref{eq:level1auto} acting on the set of transseries parameters.
However, the automorphism acting on the set of transseries parameters also involves the Stokes constants $S_{ij}$, which have the potential to change throughout the complex plane on account of the higher-order Stokes phenomenon.
This change in the Stokes constants is interpreted in \eqref{eq:HOSA+P2} as an automorphism acting on the set of Stokes constants themselves.
Our method allows for these phenomena to be determined consistently, and importantly resolves more complicated behaviour in which the value of the higher-order Stokes automorphism changes along the respective higher-order Stokes line.
This latter behaviour is associated with intersecting higher-order Stokes lines.

We note that our method has assumed that higher-order Stokes lines do not perfectly coincide with certain Stokes lines, which would result in an adjusted Stokes automorphism \citep{nemes2022dingle,shelton2024exponential}, and also that certain higher-order Stokes lines do not coincide with one another, which would result in a new form of higher-order Stokes automorphism. This is a complicated problem that depends substantially on the nature of the Stokes and higher-order Stokes lines that coincide; calculating this behaviour would be an interesting problem for future studies in this area.
The case in which asymptotic solutions are dominated by an algebraic (Poincar\'e) expansion along the real axis, as common in problems in applied mathematics and physics, is a special case of our method with $\chi_1(z)=0$.
In such cases, the effect of nonlinearity can also be resolved but here one has to be careful in selecting which transseries parameters and Stokes constants to study, since the number of both of these will be unbounded.

Our method has been demonstrated by solving for the far-field asymptotics of the swallowtail integral, as expressed as the solution to a fourth-order differential equation.
In this example we have demonstrated how the crossing of higher-order Stokes lines, and therefore a change in value of the higher-order Stokes automorphism, can emerge in a system with four exponentials in the asymptotic transseries.
This example has not been considered previously in the generality required to display this behaviour.
The result of this phenomenon displayed by the higher-order Stokes lines was seen to be a change of activity of multiple Stokes lines through the same Stokes crossing point, as seen in figure~\ref{fig:activeSL}.

\section*{Acknowledgments} 
This work was formulated during the Isaac Newton Institute programme on ``Applicable resurgent asymptotics: towards a universal theory" and the Okinawa Institute of Science and Technology visiting programme ``Exact Asymptotics: From Fluid Dynamics to Quantum Geometry".
We also acknowledge helpful discussions with Prof. G. Nemes (Harbin), Prof. C. J. Howls (Southampton), and Prof. I. Aniceto (Southampton). CJL acknowledges the support of ARC Discovery Project DP240101666. The work of SC is supported by the National Science and Technology Council of Taiwan under grants NSTC 113-2112-M-007-019 and NSTC 114-2112-M-007-015.

\vskip2pc
\bibliographystyle{jfm}
\bibliography{references}
\providecommand{\noopsort}[1]{}

\appendix

\section{Parametric Alien calculus}
\label{sec:alienapp}
We give a brief introduction to \'Ecalle's Alien calculus applied to parametric asymptotic series that arise, for example, as formal solutions to singularly perturbed linear differential equations. Let us first denote by $\psi(z,\epsilon)$ a parametric Gevrey-1 asymptotic series and write
\begin{equation}
    \psi(z,\epsilon) \sim \psi_0(z) + \epsilon \psi_1(z) + \epsilon^2 \psi_2(z) + \cdots.
\end{equation}
If we suppose further that this series is \textit{resurgent}\footnote{Away from a discrete set of `turning points' $z \in \mathbb{C}_z$.} \citep{ecalle1981iii} meaning that its Borel transform is endlessly analytically continuable\footnote{A complex power series is endlessly analytically continuable if it may be analytically continued along any path that avoids a discrete set of singular points.} then the large $n$ asymptotics of the coefficients $\psi_n(z)$ are given by the descending factorial expression
\begin{equation}\label{eq:appendixlateterms}
    \psi_n(z) \sim \sum_{\tilde{\chi}_i} \sum_{p=0}^{\infty} \tilde{\psi}^{(i)}_p(z) \frac{\Gamma(n+\alpha_j-p)}{\left(\tilde{\chi}_j(z)\right)^{n+\alpha_j-p}}.
\end{equation}
Denoting by $\mathbb{C}_1[[\epsilon]]_1$ this ring\footnote{\cite{ecalle1981iii} demonstrates that the resurgent property is closed under addition and multiplication of formal power series.} of resurgent asymptotic series, the \textit{Alien derivatives}, are a set of differential operators $\Delta: \mathbb{C}_1[[\epsilon]]_1 \to \mathbb{C}_1[[\epsilon]]_1$ associated to each component in the late term asymptotics. The set of Alien derivatives and action on $\psi(z,\epsilon)$ defined above is then
\begin{equation}
    \Delta_{\tilde{\chi}_j}: \psi(z,\epsilon) \mapsto \epsilon^{-\alpha_j}\sum_{p=0}^{\infty} \tilde{\psi}^{(j)}_p(z)\epsilon^p.
\end{equation}
Linearity 
\begin{equation}
    \Delta_{\chi}\left(a(z) \psi(z,\epsilon) + b(z)\phi(z,\epsilon)\right) = a(z) \Delta_{\chi}\psi(z,\epsilon) + b(z) \Delta_{\chi}\phi(z,\epsilon),
\end{equation}
and the Liebniz property
\begin{equation}
    \Delta_{\chi}(\psi(z,\epsilon)\phi(z,\epsilon)) = \psi(z,\epsilon) \Delta_{\chi}\phi(z,\epsilon) + \phi(z,\epsilon) \Delta_{\chi}\psi(z,\epsilon),
\end{equation}
follow directly from the definition and justify the nomenclature Alien \emph{derivative}. Finally, we define the so-called \textit{dotted} Alien derivative by $\dot{\Delta}_{\chi} = \mathrm{e}^{-\chi(z)/\epsilon} \Delta_{\chi}$. 

An immediate but key property of the dotted Alien derivative is that it commutes with the ordinary (singularly-perturbed) derivative, namely
\begin{equation}\label{eq:Aliencommute}
    \dot{\Delta}_{\chi}\left(\epsilon \frac{\mathrm{d\psi}}{\mathrm{d}z}\right) = \epsilon \frac{\mathrm{d}}{\mathrm{d}z}\left(\dot{\Delta}_{\chi}\psi\right).
\end{equation}
The action of the Alien derivative and the dotted Alien derivative extend naturally\footnote{That is, exponential terms are declared to have zero Alien derivative.} to act linearly on formal transseries. Let us now consider such a finite transseries
\begin{equation}
    \Psi(z,\epsilon;\sigma_1,\ldots,\sigma_N) := \sum_{i=1}^N \sigma_i \epsilon^{-\alpha_i}\mathrm{e}^{-\chi_i(z)/\epsilon}\left(\psi^{(i)}_0(z) + \psi^{(i)}_1(z)\epsilon + \psi^{(i)}_2(z)\epsilon^2 + \cdots \right)
\end{equation}
together with a corresponding set of Alien derivatives $\Delta_{\chi_i-\chi_j}$ labelled by differences of the transseries exponentials.
Suppose further that $\Psi(z,\epsilon;\sigma_1,\ldots,\sigma_N)$ arises as the formal transseries solution to a linear ODE, namely we consider a singularly perturbed linear operator
\begin{equation}
    \mathscr{P} = \sum_{n=0}^N a_n(z)\epsilon^n \frac{\mathrm{d}^n}{\mathrm{d}z^n}
\end{equation} 
such that $\mathscr{P} \Psi = 0$. The transseries parameter derivatives
\begin{equation}
    \frac{\partial}{\partial \sigma_j}\Psi(z,\epsilon) = \epsilon^{-\alpha_j}\mathrm{e}^{-\chi_j(z)/\epsilon}\left(\psi^{(j)}_0(z) + \psi^{(j)}_1(z)\epsilon + \psi^{(j)}_2(z)\epsilon^2 + \cdots \right),\quad j=1,2,\ldots,N
\end{equation}
are $N$ linearly independent solutions, distinguished by their far-field boundary conditions, each satisfying $\mathscr{P}\left(\frac{\mathrm{d}\Psi}{\mathrm{d}\sigma_j}\right) = 0$. Crucially, the commutative property \eqref{eq:Aliencommute} ensures that the Alien derivatives also satisfy the same ODE $\mathscr{P}\left(\dot{\Delta}_{\chi_j}\Psi\right)=0$ and therefore must be linearly dependent with derivatives of the transseries with respect to the transseries parameter. This statement is encoded compactly as the \textit{second bridge equation}:
\begin{equation}
    \dot{\Delta}_{\chi_i-\chi_j} \Psi = \frac{S_{ij}}{2 \pi i} \sigma_i \frac{\partial \Psi}{\partial \sigma_j},
\end{equation}
where the constants of proportionality define the \textit{Stokes' constants} $S_{ij}$. Finally, using the bridge equation, we note that the linear Stokes automorphisms $\mathcal{S}_{i>j}: \sigma_j \to \sigma_j + S_{ij}\sigma_i$ may be expressed as the operator
\begin{equation}
    \mathcal{S}_{i>j} = \exp(2\pi i\dot{\Delta}_{\chi_i-\chi_j}).
\end{equation}
It is in this sense that (dotted) Alien derivatives generate Stokes automorphisms.

\end{document}